\documentclass[10 pt]{article} 
\usepackage{amsmath} 
\usepackage{amsthm}
\usepackage{graphicx}

\newcommand{\floor}[1]{\lfloor {#1} \rfloor}

\newtheorem{theorem}{Theorem}[section]
\newtheorem{prop}[theorem]{Proposition}
\newtheorem{lemma}[theorem]{Lemma}
\newtheorem{cor}[theorem]{Corollary}

\theoremstyle{remark} \newtheorem*{remark}{Remark}

\theoremstyle{definition} 

\title{Fractal Sequences and Restricted Nim}
\author{Lionel Levine\footnote{Supported by a National Science Foundation 
Graduate Research Fellowship.} \\ Department of Mathematics \\ University 
of California \\ Berkeley, CA, 94720}
\date{}

\DeclareSymbolFont{AMSb}{U}{msb}{m}{n}
\DeclareMathSymbol{\N}{\mathbin}{AMSb}{"4E}
\DeclareMathSymbol{\Z}{\mathbin}{AMSb}{"5A}
\DeclareMathSymbol{\R}{\mathbin}{AMSb}{"52}
\DeclareMathSymbol{\C}{\mathbin}{AMSb}{"43}

\begin{document}

\maketitle

\begin{abstract} 
The {\it Grundy number} of an impartial game $G$ is the size of the 
unique Nim heap equal to $G$.  We introduce a new variant of Nim, 
{\it Restricted Nim}, which restricts the number of stones a player 
may remove from a heap in terms of the size of the heap.  Certain 
classes of Restricted Nim are found to produce sequences of Grundy 
numbers with a self-similar fractal structure.  Extending work of C. 
Kimberling, we obtain new characterizations of these ``fractal 
sequences'' and give a bijection between these sequences and certain 
upper-triangular arrays.  As a special case we obtain the game of {\it 
Serial Nim}, in which the Nim heaps are ordered from left to right, and 
players can move only in the leftmost nonempty heap.  
\end{abstract}

\section{Introduction}
The classic game of Nim, first studied by C. Bouton \cite{Bouton}, is
played with piles of stones.  On her turn, a player can remove any
number of stones from any one pile.  The winner is the player to take 
the last stone.  Many variants of Nim have been studied; see 
chapters 14--15 of \cite[vol. 3]{WW} as well as \cite{Endnim, AN, ONAG, 
HR, Schwenk, SS, Zieve}.  In {\it Restricted Nim}, we place an upper or 
lower bound on the number of stones that can be removed in terms of the 
size of the pile.  For example, suppose the players are permitted to 
remove any number of stones strictly smaller than half the size of the 
pile.  Then a pile of size $2^n$ is a win for the second player: no 
matter how the first player moves, the second player can respond by 
reducing the size of the pile to $2^{n-1}$; when just two stones remain, 
the first player is unable to move and loses.  Likewise, the first player 
can win from any pile whose size is {\it not\/} a power of two by 
reducing the size to a power of two.
 
In general, we may require that no more than $f(m)$ stones be 
removed from a pile of size $m$; here $f$ may be any sequence of 
nonnegative integers satisfying $f(m) \leq m$.  This is the game of 
{\it Maximum Nim}.  Since the sequence $f$ specifies the rules of the 
game, we will often refer to $f$ as the {\it rule sequence}, or simply 
the {\it rule}.
 
Maximum Nim is an example of an {\it impartial game}.  By the 
{\it Sprague-Grundy theory of impartial games} \cite{WW, ONAG, 
Grundy, Sprague} any impartial game $G$ is equal to 
a Nim heap of size $g$ for some $g$.  The integer $g$ is unique and is 
called the {\it Grundy number} of $G$. (For an explanation of 
impartial games, Sprague-Grundy theory and the notion of equality 
of games, we refer the reader to the first volume of \cite{WW}.)

For each $n$, the game of Maximum Nim with rule $f$ on a pile 
of size $n$ has a Grundy number $g_n$.  The sequence $(g_n)_{n \geq 
0}$ will be called the {\it Grundy sequence} for Maximum Nim with rule 
$f$.  By the Sprague-Grundy theory, the sequence $g_n$ satisfies the 
recurrence
	\begin{equation}
	\label{grundy}
	g_n = \text{mex} \{ g_{n-i} \}_{i=1}^{f(n)}, 
	\end{equation} 
in which mex $S$ denotes the {\it minimal excludant\/} of the set 
$S$, the smallest nonnegative integer not in $S$.  

Returning to our example in which the number of stones taken must 
be strictly smaller than half the size of the pile, the recurrence 
(\ref{grundy}) with $f(n) = \floor{\frac{n-1}{2}}$ gives the 
sequence $g_n$, starting from $n=1$, as
	\begin{equation}
	\label{sequence}	
	0, {\bf 0}, 1, {\bf 0}, 2, {\bf 1}, 3, {\bf 0}, 4, {\bf 2}, 5, 
	{\bf 1}, 6, {\bf 3}, 7, {\bf 0}, 8, {\bf 4}, 9, {\bf 2}, 10, 
	\ldots
	\end{equation}
The odd-indexed terms are just the nonnegative integers in order, 
while the even-indexed terms, shown in {\bf bold}, form a copy of 
the original sequence!  This fractal-type property is a consequence 
of Theorem \ref{main}.  Note that the zeros in the sequence occur at 
positions indexed by the powers of two; these are precisely the pile 
sizes resulting in a second player win.  

The Grundy sequence (\ref{sequence}) is an example of a 
``divide-and-conquer sequence'' \cite{Dumas}.  It appeared in 
\cite{CardSorting} in the solution to a card sorting problem.  More 
generally, it is an example of the ``fractal sequences'' studied by 
Kimberling \cite{Kimberling2, Kimberling3}.  In fact, as shown in 
Proposition \ref{fractalequivdef}, all of Kimberling's fractal 
sequences can be obtained as sequences of Grundy numbers for games of 
Maximum Nim.

An explicit formula for the Grundy sequence (\ref{sequence}) is given 
by ``truncating at the last binary one:'' if $n$ is written 
in binary as
	$$ n = 2^a + \ldots + 2^y + 2^z $$
with $a > \ldots > y > z \geq 0$, then
	\begin{equation}
	\label{chopoff}
	g_n = 2^{a-z-1} + \ldots + 2^{y-z-1}.
	\end{equation}

By placing a lower bound, rather than an upper bound, on the 
number of stones that may be taken in a turn, we obtain another 
variant, {\it Minimum Nim}.  Although the analysis of Minimum Nim is 
significantly easier than that of Maximum Nim, there is a curious 
relationship between the two games.  For example, if a move consists of 
taking {\it at least\/} half the stones in a pile, any move will reduce 
by at least one the number of binary digits in the size of the pile.  
The values of the digits ($0$ or $1$) may change, but the number of 
digits is always reduced.  

We may think of this game as played with piles of red and blue beads: a 
move consists of removing any number of beads from any one pile, and in 
addition changing the colors of any number of beads remaining in that 
pile.  Of course, color has no effect on this game, which is just Nim.  
Playing Minimum Nim with this rule on a pile of size $n$ is thus 
equivalent to playing ordinary Nim on the binary digits of $n$.  In 
other words, the Grundy number $h_n$ for a pile of size $n$ is just 
$\floor{\log_2 n} + 1$.  

Notice that in our example, $n$ can be uniquely recovered from its pair 
of Grundy numbers $(g_n, h_n)$ for Maximum and Minimum Nim: by 
comparing $h_n$ with the number of binary digits in $g_n$, we can 
determine how many final zeros were deleted when using (\ref{chopoff}) 
to pass from $n$ to $g_n$.  To recover $n$, simply write $g_n$ in 
binary and append a final one followed by the appropriate number of 
zeros.  Theorem \ref{bijection2} generalizes this observation.

\section{Maximum Nim}

When the rule sequence $f$ is weakly increasing, the corresponding Grundy 
sequence $g_n$ for Maximum Nim exhibits a self-similar fractal 
structure.  Sequences $f$ satisfying
	\begin{equation}
	\label{regular}
	0 \leq f(n)-f(n-1) \leq 1
	\end{equation}
play a special role in the analysis and will be called {\it regular}.  
The following lemma converts the recurrence
	\begin{equation}
	\label{grundyagain}
	g_n = \text{mex} \{g_{n-i}\}_{i=1}^{f(n)}
	\end{equation}
into a more explicit recurrence (\ref{fractal}).

\begin{lemma}
\label{fractalrec}
If $f$ is a regular sequence, the Grundy sequence $(g_n)_{n \geq 0}$ 
for Maximum Nim with rule $f$ satisfies
	\begin{equation}
	\label{fractal}
	 g_n = \begin{cases}
               f(n) & \text{{\em if }} f(n)>f(n-1); \\
               g_{n-f(n)-1}, & \text{{\em if }} f(n)=f(n-1).
             \end{cases}
	\end{equation}
\end{lemma}

\begin{proof}
Fix $0 \leq j \leq n$.  By regularity, $f(n) \leq j+f(n-j)$, so $g_{n-j} = 
\text{mex} \{g_{n-i}\}_{i=j+1}^{j+f(n-j)}$ is distinct from $g_{n-j-1}, 
\ldots, g_{n-f(n)}$.  Thus for any $n$ the terms $g_n$, $g_{n-1}$, 
$\ldots$, $g_{n-f(n)}$ are distinct.

If $f(n)>f(n-1)$, then for any $0 < j \leq n$, by (\ref{grundyagain}) the 
term $g_{n-j}$ is the mex of a set of size strictly smaller than $f(n)$, 
hence $g_{n-j} < f(n)$.  Since $g_{n-1}, \dots, g_{n-f(n)}$ are distinct 
and $<f(n)$, they must be $0,1, \dots, f(n)-1$ in some order.  Thus $g_n = 
\text{mex} \{0,1, \dots, f(n)-1 \} = f(n)$, completing the proof in the 
first case.

Now suppose $f(n)=f(n-1)$.  Since $g_{n-1}, \dots, g_{n-1-f(n-1)}$ 
are distinct and $\leq f(n)$, they are $0,1,\dots, f(n)$ in some 
order, so 
	\begin{equation*}
	g_{n-f(n)-1} = g_{n-1-f(n-1)} = \text{mex} 
		\{g_{n-i}\}_{i=1}^{f(n)} = g_n. \qed
	\end{equation*}
\renewcommand{\qed}{} \end{proof}
Following \cite{Kimberling3}, we denote by $\Lambda(g)$ the subsequence 
of $g$ obtained by deleting, for each integer $i \geq 0$, the first 
term equal to $i$.  As the following theorem shows, the Grundy 
sequences for Maximum Nim are ``self-similar'' in the sense that they 
satisfy $\Lambda(g) = g$. 

\begin{theorem}
\label{main}
Let $f$ be a regular sequence, and let $(g_n)_{n \geq 0}$ be the Grundy 
sequence for Maximum Nim with rule $f$.  Then $\Lambda(g) = g$. 
\end{theorem}

\begin{proof}
By Lemma \ref{fractalrec}, $\Lambda(g)$ consists of precisely those 
terms $g_n$ for which $f(n) = f(n-1)$.  Since $f$ is regular, it follows 
that all but $f(n)+1$ of the terms $g_0, g_1, \ldots, g_n$ lie in the 
subsequence $\Lambda(g)$.  Thus if $f(n)=f(n-1)$, we have by Lemma 
\ref{fractalrec}
	\begin{equation}
	\label{lambda}
	\Lambda(g)_{n-f(n)-1} = g_n = g_{n-f(n)-1}.
	\end{equation}
Since $f$ is regular, as $n$ ranges through all positive integers such 
that $f(n) = f(n-1)$, the quantity $n-f(n)-1$ ranges through all 
nonnegative integers, and hence $\Lambda(g) = g$.
\end{proof}

Lemma \ref{fractalrec} and Theorem \ref{main} provide an easy algorithm 
for writing down the first $n$ terms of the Grundy sequence $g$ in time 
$O(n)$. (This is a significant improvement over the recurrence 
(\ref{grundyagain}), which requires time on the order of $\sum_{i=1}^n 
f(i)$.)  First, make a table of the values $f(0), \dots, f(n)$, marking 
those indices $n_1 < \ldots < n_k$ for which $f(n_i) > f(n_i - 1)$.  
Next, write the integers $0, 1, \ldots, k$ in positions $0, n_1, \dots, 
n_k$; this takes care of the first case in (\ref{fractal}).  Finally, 
fill in the gaps between the $n_i$ in the unique way possible so that 
the gapped sequence forms a copy of the original; this is done by 
copying earlier terms according to the second case of (\ref{fractal}).  
The example below illustrates this algorithm for the rule sequence $f(n) = 
\floor{\sqrt{n}}$.

	\begin{center}
        
\begin{tabular}{r@{~~~}r@{~~~}r@{~~~}r@{~~~}r@{~~~}r@{~~~}r@{~~~}r@{~~~}r@{~~~}r@{~~~}r@{~~}r@{~~}r@{~~}r@{~~}r@{~~}r@{~~}r@{~~}r@{~~}r@{~~}r@{~~}r@{~~}r@{~~}} 
                 $n$&0&1&2&3&4&5&6&7&8&9&10&11&12&13&14&15&16\\
              $f(n)$&0&1&1&1&2&2&2&2&2&3& 3& 3& 3& 3& 3& 3& 4\\
     case 1&0&1& & &2& & & & &3&  &  &  &  &  &  & 4\\
     case 2& & &{\bf 0}&{\bf 1}& &{\bf 0} &{\bf 1} &{\bf 2} & 
{\bf 0} & & {\bf 1}&{\bf 2} &{\bf 0} &{\bf 3} &{\bf 1} &{\bf 2}&  
        \end{tabular}
	\end{center}

Since $\floor{\sqrt{n}}$ exceeds $\floor{\sqrt{n-1}}$ precisely 
when $n$ is a perfect square, we have $n_i = i^2$.  The first case 
of (\ref{fractal}) gives $g_{i^2} = i$, and the second case is 
used to compute the remaining terms.

Our next result reduces the problem of computing Grundy numbers for a
general weakly increasing rule sequence $f$ to the case of regular $f$, 
so that Theorem \ref{main} applies.

\begin{prop}
\label{extension}
If $f$ is any weakly increasing sequence, the Grundy sequence for 
Maximum Nim with rule $f$ is the same as that with rule $f'$, where 
the regular sequence $f'$ is defined inductively by
	$$ f'(n) = \min \{f(n), 1+f'(n-1)\}. $$  
\end{prop}

\begin{proof}
Let $g_n$ and $g'_n$ be the Grundy sequences corresponding to rules 
$f$ and $f'$, and induct on $n$ to show $g_n = g'_n$.  If $f'(n)=f(n)$, 
then the inductive hypothesis, together with (\ref{grundyagain}), implies 
$g_n = g'_n$.  Otherwise, $f'(n) = 1+f'(n-1) < f(n)$.  Since $f'$ is 
regular, by Lemma \ref{fractalrec} we have $g'_n = f'(n) > g'_{n-j} = 
g_{n-j}$ for all $0 < j \leq n$, hence
	$$ g'_n \geq \text{mex} \{g_{n-j}\}_{j=1}^{f(n)}
	 	\geq \text{mex} \{g_{n-j}\}_{j=1}^{f'(n)}
		   = \text{mex} \{g'_{n-j}\}_{j=1}^{f'(n)}
	           = g'_n, $$
hence $g'_n = \text{mex} \{g_{n-j}\}_{j=1}^{f(n)} = g_n$.
\end{proof}

By way of example, consider the rule sequence $f(n) = \text{max} \{2^k 
\leq n \} -1$: players may remove any number of stones less than the 
greatest power of two not exceeding the size of the pile.  Since $f$ is 
not regular, we use Proposition \ref{extension} to pass to the regular 
sequence $f'$ before applying Theorem \ref{main}.  The following chart 
gives values for $f$, $f'$ and $g$.

	\begin{center}
	
\begin{tabular}{r@{~~~}r@{~~~}r@{~~~}r@{~~~}r@{~~~}r@{~~~}r@{~~~}r@{~~~}r@{~~~}r@{~~~}r@{~~}r@{~~}r@{~~}r@{~~}r@{~~}r@{~~}r@{~~}r@{~~}r@{~~}r@{~~}r@{~~}r@{~~}r@{~~}}	
	   $n$&0&1&2&3&4&5&6&7&8&9&10&11&12&13&14&15&16\\
	$f(n)$&0&0&1&1&3&3&3&3&7&7& 7& 7& 7& 7& 7& 7&15\\
       $f'(n)$&0&0&1&1&2&3&3&3&4&5& 6& 7& 7& 7& 7& 7& 8\\ 
	 $g_n$&0&{\bf 0}&1&{\bf 0}&2&3&{\bf 1}&{\bf 0}&4&5& 6& 7& 
	{\bf 2}& {\bf 3}& {\bf 1}& {\bf 0}& 8
	\end{tabular}
	\end{center}

If $n$ is one less than a power of two, then $g_n=0$.  Otherwise, 
writing $n$ in binary, after the inital $1$ there will be a string 
of ones, possibly empty, followed by a zero: $n = (1 1^k 0 b_1 
\dots b_j)_2$.  Now $g_n$ is obtained by deleting this string of ones
and the zero that follows it: $g_n = (1 b_1 \dots b_j)_2$.

\section{Fractal Sequences}

We now show that the Grundy sequences for Maximum Nim with a weakly 
increasing rule $f$ are precisely the ``fractal sequences'' studied by 
Kimberling \cite{Kimberling2, Kimberling3}.  Following 
\cite{Kimberling3}, we call a sequence $(g_n)_{n \geq 0}$ {\it 
infinitive} if for every integer $k \geq 0$ infinitely many terms $g_n$ 
are equal to $k$.  A {\it fractal sequence} $(g_n)_{n \geq 0}$ is 
defined in \cite{Kimberling3} as an infinitive sequence satisfying two 
additional properties: \\

(F2) If $j<k$, the first instance of $j$ in $g$ precedes the first 
instance of $k$;

(F3) The subsequence $\Lambda(g)$ of $g$ obtained by deleting the first 
instance of each integer $k$ is $g$ itself. \\

By an {\it instance} of an integer $k$ in $g$ we mean a term $g_n = k$.  
If $g$ is an infinitive squence, denote by $\hat{g}(k)$ the position of 
the first instance of $k$ in $g$.  If $g$ is fractal, the sequence 
$\hat{g}$ is increasing by property (F2).

\begin{lemma}
\label{hat}
Let $g$ and $h$ be fractal sequences.  If $\hat{g}=\hat{h}$, then $g=h$.
\end{lemma}

\begin{proof}
Induct on $n$ to show $g_n = h_n$.  If $n = \hat{g}(k)$ for some $k$, 
then $g_n = h_n = k$.  Otherwise, let $k$ be such that $\hat{g}(k) < n < 
\hat{g}(k+1)$.  By property (F3) and the inductive hypothesis,
	$$ g_n = \Lambda(g)_{n-k-1} = g_{n-k-1} = h_{n-k-1} = 
\Lambda(h)_{n-k-1} = h_n. \qed $$
\renewcommand{\qed}{} 
\end{proof}

\begin{prop}
\label{fractalequivdef}
Let $(g_n)_{n \geq 0}$ be an infinitive sequence.  The following are 
equivalent.

{\em (i)} $g$ is a fractal sequence;

{\em (ii)} $g$ is the Grundy sequence for Maximum Nim for some weakly 
increasing rule sequence $f$;

{\em (iii)} $g$ is the Grundy sequence for Maximum Nim for some regular 
rule sequence $f$.
\end{prop}

\begin{remark}
Equivalently, conditions, (ii) and (iii) may be replaced by the 
condition that $g$ satisfies the recurrence
	$$ g_n = \text{mex} \{g_{n-i}\}_{i=1}^{f(n)} $$
for a weakly increasing or regular sequence $f$, respectively.
\end{remark}

\begin{proof}
We'll show (iii) $\Rightarrow$ (ii) $\Rightarrow$ (i) $\Rightarrow$ 
(iii).  The first implication is trivial.  If $g$ is the Grundy 
sequence for Maximum Nim with rule $f$, with $f$ weakly increasing, by 
Proposition \ref{extension} it follows that $g$ is also the Grundy 
sequence for rule $f'$, which is regular.  By Theorem \ref{main}, it 
follows that $g$ is a fractal sequence.  

For the final implication, let $f(n) = \max \{g_m\}_{m \leq n}$.  By 
property (F2), we have $0 \leq f(n) - f(n-1) \leq 1$, i.e.\ $f$ is 
regular.  Let $h_n$ be the Grundy sequence for Maximum Nim with 
rule $f$.  We will show $g=h$.  By Theorem \ref{main}, $h$ is a fractal 
sequence, and by Lemma \ref{fractalrec}
	\begin{eqnarray*}
	\hat{h}(k) &=& \min \{n | f(n)=k\} \\
		   &=& \min \{n | \max \{g_m\}_{m \leq n} = k \} \\
		   &=& \min \{n | g_n = k \} \\
		   &=& \hat{g}(k).
	\end{eqnarray*}
By Lemma \ref{hat} it follows that $g=h$.
\end{proof}

Kimberling \cite{Kimberling1, Kimberling3} has given characterizations 
of fractal sequences---the notions of {\it interspersion} and {\it 
dispersion}---which on the surface have nothing to do with 
self-similarity.  These are defined in terms of an {\it associated 
array} \cite{Kimberling3} $A = A(g) = (a_{ij})_{i,j \geq 0}$ whose 
$i$-th row consists of the instances of $i$ in $g$ listed in increasing 
order.  The array $A(g)$ contains every positive integer exactly once, 
and its rows are increasing.  An array having these properties is 
called an {\it interspersion} if, in addition, its columns are 
increasing and \\

(I4) $a_{ij} < a_{kl} < a_{i,j+1}$ implies $a_{i,j+1} < a_{k,l+1} <
a_{i,j+2}$. \\

\noindent In \cite{Kimberling3} it is shown that

\begin{theorem} 
\label{Kimberling}
$g$ is a fractal sequence if and only if $A(g)$ is an interspersion. 
\end{theorem}
 
We find it illuminating to recast the definition of an interspersion in 
terms of the sequence itself, rather than its associated array.  If $M$ 
is a set of nonnegative integers and $g$ an infinitive sequence, the 
{\it restriction of $g$ to $M$}, denoted $g|M$, is the subsequence of 
$g$ formed by deleting all terms $g_n$ for which $g_n \notin M$.  In 
these terms, an {\it interspersion} is an infinitive sequence $g$ such 
that for any $i<j$ the restriction $g|\{i,j\}$ has the form
	$$ i, i, i, \ldots, i, j, i, j, i, j, \ldots; $$
after an initial segment of $i$'s, instances of $i$ and $j$ must 
alternate.

When $M$ is infinite, it is often useful to relabel the sequence $g|M$ 
so as to make it infinitive.  If $M = \{m_0, m_1, \ldots \}$ with $m_0 
< m_1 < \ldots$, the {\it relabeling} of $g|M$ is the sequence obtained 
by replacing each instance of $m_i$ with $i$.

Our next result characterizes the restrictions of an interspersion.  
Taking $M$ to be the set of positive integers, we obtain as a special 
case Theorem 5 of \cite{Kimberling3}.

\begin{prop}
Let $g$ be an interspersion, and let $M$ be a set of nonnegative 
integers.

	{\em (i)} If $M$ is finite, then $g|M$ is eventually periodic 
with period $\#M$.

	{\em (ii)} If $M$ is infinite, the relabeling of $g|M$ is an 
interspersion.
\end{prop}

\begin{proof}
(i) Let $m=\#M$, and fix $i \in M$.  With finitely many exceptions, 
between consecutive instances of $i$ in $g|M$ there is exactly one 
instance of each $j \in M-\{i\}$.  Thus for sufficiently large $n$ the 
$m$ terms
	$$ (g|M)_n, (g|M)_{n+1}, \ldots, (g|M)_{n+m-1} $$
are a permutation of $M$.  In particular, both $(g|M)_n$ and 
$(g|M)_{n+m}$ are equal to the unique element $j \in M$ not contained 
in $\{(g|M)_{n+i}\}_{i=1}^{m-1}$, so $g|M$ is eventually periodic mod 
$m$.

(ii) Write $M = \{m_0, m_1, \ldots \}$ with $0 \leq m_0 < m_1 < 
\ldots$.  For $i<j$, since the restriction $g|\{m_i, m_j\}$ has the 
form
	$$ m_i, m_i, \ldots, m_i, m_j, m_i, m_j, \ldots, $$
the restriction of the relabeling of $g|M$ to $\{i,j\}$ has the form
	$$ i,i, \ldots, i, j, i, j, \ldots, $$
so the relabeling of $g|M$ is an interspersion.
\end{proof}

If $g$ is an interspersion, the restriction $g|\{i,j\}$ is determined 
by the number $s_{ij}$ of instances of $i$ in $g$ preceding the first 
instance of $j$.  (If $i=0$, we do not count the instance $g_0 = 0$.)  
The array $S(g) := (s_{ij})_{i,j \geq 0}$ is strictly upper-triangular 
and satisfies
	\begin{equation}
	\label{subadditive}
	s_{ij} + s_{jk} - 1 \leq s_{ik} \leq s_{ij} + s_{jk}.
	\end{equation}
Equality holds on the left or the right side of (\ref{subadditive}) 
accordingly as the restriction $g|\{i,j,k\}$ has the form
     $$ i, i, \ldots, i, j, i, j, \ldots, i, j, k, i, j, k, i, \ldots $$
or
     $$ i, i, \ldots, i, j, i, j, \ldots, i, j, i, k, j, i, k, \ldots.$$
An upper-triangular array satisfying (\ref{subadditive}) will be called 
a {\it subadditive triangle.}  

For example, the array
	$$\begin{array}{ccccccccccc}
 	2& 3& 3& 4& 4& 4& 4& 5& 5& 5& \ldots \\
  	 & 1& 2& 2& 2& 3& 3& 3& 3& 3& \ldots \\
  	 &  & 1& 1& 2& 2& 2& 2& 2& 3& \ldots \\
  	 &  &  & 1& 1& 1& 2& 2& 2& 2& \ldots \\
  	 &  &  &  & 1& 1& 1& 1& 2& 2& \ldots \\
  	 &  &  &  &  & 1& 1& 1& 1& 1& \ldots \\
  	 &  &  &  &  &  & 1& 1& 1& 1& \ldots \\
  	 &  &  &  &  &  &  & 1& 1& 1& \ldots \\
 	 &  &  &  &  &  &  &  & 1& 1& \ldots \\
  	 &  &  &  &  &  &  &  &  & 1& \ldots
	\end{array}$$
is the subadditive triangle associated to the Grundy sequence
	$$ 0, {\bf 0}, 1, {\bf 0}, 2, {\bf 1}, 3, {\bf 0}, 4, {\bf 2}, 
	5, {\bf 1}, 6, {\bf 3}, 7, {\bf 0}, 8, {\bf 4}, 9, {\bf 2}, 10, 
	\ldots $$
for Maximum Nim with rule $f(n) = \floor{\frac{n-1}{2}}$.

In Theorem \ref{bijection1} we show that the correspondence between 
fractal sequences and subadditive triangles is a bijection.  
\begin{lemma}
\label{columnsums}
A subadditive triangle $(s_{ij})_{i,j \geq 0}$ is determined by its 
column sums $c_j = \sum_{i=0}^{j-1} s_{ij}$.
\end{lemma}

\begin{proof}
For $i<j<k$ write
	$$ \varepsilon_{ijk} = s_{ij} + s_{jk} - s_{ik}. $$
By subadditivity (\ref{subadditive}), each $\varepsilon_{ijk}$ is 
either $0$ or $1$.  For fixed $i$, we induct on $j$ to show that the 
column sums $c_i, c_{i+1}, \dots, c_j$ determine the entry $s_{ij}$.  
We have 
	\begin{eqnarray}
  c_j &=& \sum_{p=0}^{j-1} s_{pj} \nonumber \\
      &=& \sum_{p=0}^{i-1} (s_{pi} + s_{ij} - \varepsilon_{pij})
	+ s_{ij}
	+ \sum_{p=i+1}^{j-1} (s_{ij} - s_{ip} + \varepsilon_{ipj})
		\nonumber \\
      &=& c_i + js_{ij} - \sum_{p=i+1}^{j-1} s_{ip} + \varepsilon,
	\label{breakdown}
	\end{eqnarray}
where the error term
	$$ \varepsilon = \sum_{p=i+1}^{j-1} \varepsilon_{ipj} 
		- \sum_{p=0}^{i-1} \varepsilon_{pij} $$
is bounded by
	\begin{equation}
	\label{errorbounds}
	-i \leq \varepsilon \leq j-1-i.
	\end{equation}
By the induction hypothesis, the sum 
	$$ \Sigma := \sum_{p=i+1}^{j-1} s_{ip} $$
appearing in (\ref{breakdown}) is determined by the column sums $c_i, 
c_{i+1}, \dots, c_{j-1}$.  Solving (\ref{breakdown}) for $s_{ij}$ we 
obtain
	\begin{equation}
	\label{fraction}
   s_{ij} = \frac{1}{j} \left[ c_j - c_i + \Sigma - \varepsilon \right]
	\end{equation}
in which every term on the right hand side, except the error term 
$\varepsilon$, is determined by the column sums.  By the bounds 
(\ref{errorbounds}), there is a unique value of $\varepsilon$ making 
the right hand side an integer, and hence $s_{ij}$ is determined by the 
column sums $c_i, c_{i+1}, \ldots, c_j$.
\end{proof}

\begin{theorem}
\label{bijection1}
The map $g \mapsto S(g)$ is a bijection between fractal sequences and 
subadditive triangles.
\end{theorem}

\begin{proof}
Given a fractal sequence $g$, write $s_{ij}$ for the typical entry 
of $S(g)$.  For fixed $j$, the column sum
	 $$ c_j = \sum_{i=0}^{j-1} s_{ij} $$
counts each term preceding the first instance of $j$ in $g$ exactly 
once.  Thus $\hat{g}(j) = 1 + c_j$.  By Lemma \ref{hat}, the sequence 
$\hat{g}$ determines $g$, so the map $g \mapsto S(g)$ is 1--1.

To show that the map is onto, given a subadditive triangle $S = 
\{s_{ij}\}$, let $g$ be the unique fractal sequence satisfying
	$$ \hat{g}(j) = 1 + \sum_{i=0}^{j-1} s_{ij}. $$  
Then $S(g)$ and $S$ have the same column sums $c_j$.  By Lemma 
\ref{columnsums}, it follows that $S=S(g)$.
\end{proof}

\section{Minimum Nim}
In the game of Minimum Nim with rule $f$, a move consists of removing 
{\it strictly} more than $f(m)$ stones from a pile of size $m$.  In 
Maximum Nim, on the other hand, taking exactly $f(m)$ stones is 
permitted.  The effect of this convention is to simplify the statements 
of Proposition \ref{minimax} and Theorem \ref{bijection2}, which 
describe the relationship between Minimum and Maximum Nim.  The Grundy 
sequence $(h_n)_{n \geq 0}$ for Minimum Nim obeys the recurrence
	\begin{equation}
	\label{minimex}
	h_n = \text{mex} \{h_i\}_{i=0}^{n-f(n)-1}.
	\end{equation}
If $f$ is a regular sequence, the sequence $(n-f(n))_{n \geq 0}$ is 
also regular.  To avoid trivialities that arise when this sequence is 
eventually constant, we require that
	\begin{equation}
	\label{unbounded}
	n - f(n) \rightarrow \infty
	\end{equation}
as $n \rightarrow \infty$.  Proposition \ref{mininim} solves the game of 
Minimum Nim with rule $f$ in the case that $f$ is a regular sequence 
satisfying (\ref{unbounded}).

Recall the notation $\hat{h}(n) = \min \{k: h_k=n\}$.  If $h$ is 
regular, the sequence $\hat{h}$ determines $h$.

\begin{prop}
\label{mininim}
Let $f$ be a regular sequence satisfying {\em (\ref{unbounded})}, and 
let $(h_n)_{n \geq 0}$ be the Grundy sequence for Minimum Nim with rule 
$f$.  Then $h$ is a regular sequence, $\hat{h}(0)=0$ and
	$$ \hat{h}(n) = q(\hat{h}(n-1)), $$
where
	\begin{equation}
	\label{q}
	q(k) = \min \{j : j-f(j) > k \}.
	\end{equation}
\end{prop}

\begin{proof}
Let $S_n = \{h_0, h_1, \ldots, h_{n-f(n)-1}\}$.  Since $f$ is regular, 
$S_{n-1} \subset S_n$ and $S_n$ contains at most one element not in 
$S_{n-1}$.  By (\ref{minimex}), $h_n = \text{mex}~ S_n$ and hence
	$$ h_{n-1} \leq h_n \leq 1+h_{n-1},$$
i.e.\ $h$ is regular.  Since $h_0=0$ we have $\hat{h}(0)=0$ and  
	\begin{eqnarray*}
 	\hat{h}(n) &=& \min \{k : \text{mex} 
		\{h_i\}_{i=0}^{k-f(k)-1} = n \} \\
		   &=& \min \{k : h_{k-f(k)-1} = n-1 \} \\
		   &=& \min \{k : k-f(k)-1 \geq \hat{h}(n-1) \} \\
		   &=& q(\hat{h}(n-1)). \qed
	\end{eqnarray*}
\renewcommand{\qed}{} \end{proof}

For example, if $f(n) = \floor{\frac{n-1}{2}}$ then $q(k) = 2k$, and 
Proposition \ref{mininim} gives the corresponding Grundy sequence $h$ 
for Minimum Nim as
	$$ 0,1,2,2,3,3,3,3,4,4,4,4,4,4,4,4,5,\ldots; $$
as we remarked in the introduction, its $n$-th term is $\floor{\log_2 n} + 
1$.

The following lemma, which explains the importance of the function $q$, 
is closely related to the fact that fractal sequences are also {\it 
dispersions} \cite{Kimberling1}.  

\begin{lemma}
\label{dispersion}
Let $f$ be a regular sequence satisfying {\em (\ref{unbounded})}, and let 
$(g_n)_{n \geq 0}$ be the Grundy sequence for the corresponding game of 
Maximum Nim.  With $q$ as in {\em (\ref{q})}, we have $g_{q(n)} = g_n$.
\end{lemma}

\begin{proof}
Since $(n-f(n))_{n \geq 0}$ is a regular sequence, by (\ref{q}) we have 
$q(n) - f(q(n)) = n+1$ and $q(n)-1 - f(q(n)-1) = n$, hence $f(q(n)) = 
f(q(n)-1)$.  
By Lemma \ref{fractalrec}, it follows that
	$$ g_{q(n)} = g_{q(n)-f(q(n))-1} = g_n. \qed $$
\renewcommand{\qed}{}
\end{proof}

Our next proposition relates the Grundy sequences for Minimum and Maximum 
Nim.  We write $q^0(n) = n$, $q^i(n) = q(q^{i-1}(n))$.

\begin{prop}
\label{minimax}
Let $f$ be a regular sequence satisfying {\em (\ref{unbounded})}, and 
let $(g_n)_{n \geq 0}$ and $(h_n)_{n \geq 0}$ be the Grundy 
sequences for Maximum and Minimum Nim with rule $f$.  Then
	$$ h_n = \# \{0 < k \leq n : g_k = 0 \}. $$
\end{prop}

\begin{proof}
Let $z_0=0$, and let $z_i$ be the first instance of zero in $g$ following 
$z_{i-1}$.  We'll show $z_i = q^i(0)$, where $q$ is given by (\ref{q}).  
By Lemma \ref{dispersion}, we have $g_{q^i(0)} = 0$ for all $i$.  
Conversely, suppose $g_m = 0$ for some $m>0$.  By Lemma \ref{fractalrec}, 
$f(m)=f(m-1)$, hence $m-f(m) > m-1-f(m-1)$ and $m = q(m-f(m)-1)$ by 
(\ref{q}).  Then $g_{m-f(m)-1} = g_m = 0$ by Lemma \ref{dispersion}, and 
by induction it follows that $m = q^i(0)$ for some $i$.  Proposition 
\ref{mininim} now implies that $z_i = \hat{h}(i)$.  Since $h$ is regular,
	$$ h_n = \max \{i : \hat{h}(i) \leq n \}
	       = \max \{i : z_i \leq n \}
	       = \# \{0 < k \leq n : g_k = 0 \}. \qed $$
\renewcommand{\qed}{} \end{proof}

Our next result shows that $n$ can be uniquely recovered from the pair 
$(g_n, h_n)$.

\begin{theorem}
\label{bijection2}
Let $f$ be a regular sequence satisfying {\em (\ref{unbounded})}, and 
let $(g_n)_{n \geq 0}$ and $(h_n)_{n \geq 0}$ be the Grundy 
sequences for Maximum and Minimum Nim with rule $f$.  Let 
$\{s_{ij}\}_{i,j \geq 0}$ be the subadditive triangle associated to the 
sequence $g$.  The map $n \mapsto (g_n, h_n)$ is a bijection between 
nonnegative integers and pairs $(i,j)$ of nonnegative integers 
satisfying $j \geq s_{0i}$.  
\end{theorem}

\begin{proof} $g$ is a fractal sequence by Proposition 
\ref{fractalequivdef}, and hence an interspersion by Theorem 
\ref{Kimberling}.  Thus if $g_m = g_n$ for some $m<n$, there is some term 
$g_i=0$ with $m < i \leq n$.  By Proposition \ref{minimax} it follows that 
$h_m < h_n$, hence the map $n \mapsto (g_n, h_n)$ is 1--1.  

Since $g$ is an interspersion, instances of $0$ and $i$ in $g$ alternate 
after the first instance of $i$, so by Proposition \ref{minimax}, for 
every $j \geq s_{0i}(g)$ there is an index $n$ such that $g_n = i$ and 
$h_n = j$.
\end{proof}

\begin{cor}
\label{array}
The array $A=(a_{ij})_{i,j \geq 0}$ whose entry $a_{ij}$ is the unique 
integer $n$ such that $i=g_n$, $j=h_n-s_{0i}$ is an interspersion. 
\end{cor}

\begin{proof}
The entry $a_{ij}$ of $A$ is the position of the $j$-th instance of $i$ 
in $g$; i.e.\ $A$ is the associated array of $g$.  By Theorem 
\ref{Kimberling}, since $g$ is a fractal sequence, $A$ is an 
interspersion.
\end{proof}

The array $A'=\{a'_{ij}\}_{j \geq s_{0i}}$ shown below is the inverse 
to the map $n \mapsto (g_n, h_n)$ for the rule sequence $f(n) = 
\floor{\frac{n-1}{2}}$.  The entry $a'_{ij}$ is the unique integer 
$n$ for which $i=g_n$, $j=h_n$.  The blank spaces in the lower left 
correspond to pairs $(i,j)$ satisfying $j<s_{0i}$, for which no such 
$n$ exists.
	$$\begin{array}{cccccccc}
	0&  1&  2&  4&  8& 16& 32& \ldots \\
   	 &   &  3&  6& 12& 24& 48& \ldots \\
   	 &   &   &  5& 10& 20& 40& \ldots \\
   	 &   &   &  7& 14& 28& 56& \ldots \\
   	 &   &   &   &  9& 18& 36& \ldots \\
	 &   &   &   & 11& 22& 44& \ldots \\
	 &   &   &   & 13& 26& 52& \ldots \\
	 &   &   &   & 15& 30& 60& \ldots \\
	 &   &   &   &   & 17& 34& \ldots \\
   	 &   &   &   &   & 19& 38& \ldots \\
   	 &   &   &   &   &\vdots&\vdots
	\end{array}$$ 
If the rows of $A'$ are left-justified, by Corollary \ref{array} the 
resulting array $A$ is an interspersion.

\section{Serial Nim}

In general, it seems difficult to describe the behavior of the Grundy 
sequences for Maximum and Minimum Nim when the rule sequence $f$ is not 
weakly increasing.  Certain special cases are of interest, however. In the 
game of {\it Serial Nim}, heaps are arranged in a row from left to 
right, and players can remove stones only from the leftmost nonempty 
heap.  If the heaps have sizes $a_1, \ldots, a_k$, we denote the Grundy 
number of the resulting game by $[a_1, \ldots, a_k]$.  This bracket is 
``right-associative'' in the sense that $[a_1, \ldots, a_k] = 
[a_1,[a_2, \ldots, a_k]]$. (However, it is {\it not\/} 
left-associative!)  If $f$ is a rule sequence of the form
	$$ 1, 2, \ldots, a_1, 1, 2, \ldots, a_2, \ldots, $$
then a single heap of size $n$ in the corresponding game of Maximum Nim 
is equivalent to a row of heaps of sizes $n-\sum_{i=1}^k a_i, a_k, 
a_{k-1}, \ldots, a_1$ in Serial Nim, where $k$ is such that $\sum_{i=1}^k 
a_i < n \leq \sum_{i=1}^{k+1} a_i$.
 
Consider the case of two heaps of sizes $a$, $b$.  Since $[0,b] = 
b$ and
	$$[a,b] = \text{mex} \{[i,b]\}_{0 \leq i < a},$$
by induction on $a$ the sequence $([a,b])_{a \geq 0}$ has the form
	$$ b, 0, 1, \ldots, b-1, b+1, b+2, \ldots; $$
in other words, for $a>0$ the bracket $[a,b]$ is $a-1$ or $a$ 
accordingly as $a \leq b$ or $a>b$.

Our next result treats the general case of $k$ heaps.  As with two 
heaps, the Grundy number of the game is always equal either to the 
size $a_1$ of the first heap or to $a_1 -1$. Moreover if the heap 
in position $m$ is the leftmost heap whose size differs from the 
first, then the Grundy number depends only on the parity of $m$ and the 
relative size of $a_m$ and $a_1$.  In this respect, Serial Nim behaves 
like a simplified version of the game ``End-nim'' studied by Albert and 
Nowakowski \cite{Endnim}, in which players may remove stones from either 
the leftmost or the rightmost nonempty heap.  Although the End-nim 
positions of Grundy number zero were classified in \cite{Endnim}, in 
general its Grundy numbers seem to behave erratically.  By contrast, the 
following result completely characterizes the Grundy numbers for Serial 
Nim.

\begin{prop}
Let $a_1, \dots, a_k$ be positive integers, and set $a_{k+1}=0$.  Let
$m = \min \{j|a_j \neq a_1\}$.  If $m$ is odd and $a_m<a_1$, or $m$
is even and $a_m>a_1$, then $[a_1, \dots, a_k] = a_1 -1$; otherwise
$[a_1, \dots, a_k] = a_1$.
\end{prop}

\begin{proof}
Induct on $k$.  The base case $k=2$ is discussed above.  Write $a 
= [a_1, \dots, a_k] = [a_1, b]$, where $b = [a_2, \dots, a_k]$.  
By the inductive hypothesis, $a = [a_1, b] = a_1 -1$ or $a_1$ 
accordingly as $b \geq a_1$ or $b \leq a_1-1$.  If $m$ is odd, 
then $a_2 = a_1$ and by the inductive hypothesis $b = a_2 -1$ or 
$a_2$ accordingly as $a_m>a_2$ or $a_m<a_2$, i.e.\ $a = a_1 -1$ or 
$a_1$ accordingly as $a_m<a_1$ or $a_m>a_1$.

Suppose now that $m$ is even.  If $m=2$, then either $a_2<a_1$, in 
which case $b \leq a_2 \leq a_1 -1$, so $a = a_1$; or $a_2>a_1$, in 
which case $b \geq a_2-1 \geq a_1$, hence $a=a_1 -1$.  If $m>2$, then 
$a_2 = a_1$ and $b = a_2 -1$ or $a_2$ accordingly as $a_m<a_2$ or 
$a_m>a_2$, i.e.\ $a = a_1 - 1$ or $a_1$ accordingly as $a_m>a_1$ or 
$a_m<a_1$.
\end{proof}

A closely related game is ``Smallest Nim,'' \cite[v.\ 3]{WW} in 
which players may take stones only from the heap (or one of the heaps) of 
smallest size.  Smallest Nim is the special case of Serial Nim in which 
the piles are arranged in nondecreasing order of size.  Further Nim 
variants in which moves are permitted to occur in only one pile are 
studied in \cite{AN}.

\section*{Acknowledgments} The author would like to thank Prof. Elwyn 
Berlekamp for helpful suggestions regarding content and exposition.

\end{document}